\documentclass[12pt,twoside]{amsart}

\usepackage{amsfonts,amsmath,amssymb,latexsym,graphicx,amscd}

\newtheorem{thm}{Theorem}[section]
\newtheorem{prop}[thm]{Proposition}
\newtheorem{lem}[thm]{Lemma}
\newtheorem{cor}[thm]{Corollary}
\newtheorem{claim}[thm]{Claim}
\newtheorem{prob}[thm]{Problem}
\newtheorem{ques}[thm]{Question}
\numberwithin{equation}{section}

\theoremstyle{definition}
\newtheorem{defn}[thm]{Definition}
\newtheorem{conj}[thm]{Conjecture}
\newtheorem{exmp}[thm]{Example}
\newtheorem{rem}[thm]{Remark}

\newcommand{\Hom}{\operatorname{Hom}}
\newcommand{\sHom}{\mathcal Hom}
\newcommand{\sExt}{\mathcal Ext}
\newcommand{\Ext}{\operatorname{Ext}}
\newcommand{\im}{\operatorname{im}}
\newcommand{\Spec}{\operatorname{Spec}}
\newcommand{\Pic}{\operatorname{Pic}}

\newcommand{\Coker}{\operatorname{Coker}}
\newcommand{\Hilb}{\operatorname{Hilb}}
\newcommand{\Bs}{\operatorname{Bs}}
\newcommand{\Rat}{\operatorname{Rat}}
\renewcommand{\div}{\operatorname{div}}

\newcommand{\res}{res}

\newcommand{\red}{\operatorname{red}}
\newcommand{\id}{\operatorname{id}}

\newcommand{\mapright}[1]{%
\smash{\mathop{%
\hbox to 1cm{\rightarrowfill}}\limits^{#1}}}
\newcommand{\mapleft}[1]{%
\smash{\mathop{%
\hbox to 1cm{\leftarrowfill}}\limits_{#1}}}
\newcommand{\mapdown}[1]{\Big\downarrow
\rlap{$\vcenter{\hbox{$\scriptstyle#1\,$}}$ }}

\title[Obstructions and non-reduced components]
{Obstructions to deforming space curves and 
non-reduced components of the Hilbert scheme}
\author{Hirokazu Nasu}
\address{Research Institute for Mathematical Sciences, 
         Kyoto University, 
         Kyoto, 606-8502, Japan}
\email{nasu@kurims.kyoto-u.ac.jp}

\begin{document}

\begin{abstract}
Let $H_{\mathbb P^3}^S$ denote the Hilbert scheme 
of smooth connected curves in $\mathbb P^3$.
We consider maximal irreducible closed subsets 
$W \subset H_{\mathbb P^3}^S$ whose general member $C$ 
is contained in a smooth cubic surface
and investigate the conditions for $W$ to be a component of 
$(H_{\mathbb P^3}^S)_{\red}$.
We especially study the case where the dimension of the
tangent space of $H_{\mathbb P^3}^S$ at $[C]$
is greater than $\dim W$ $(\ge 4\deg(C))$ by one.
We compute obstructions to deforming $C$ in $\mathbb P^3$
and prove that for every $W$ in this case, 
$H_{\mathbb P^3}^S$ is non-reduced along $W$
and $W$ is a component of $(H_{\mathbb P^3}^S)_{\red}$.
\end{abstract}

\maketitle

\section{Introduction}

Mumford \cite{Mumford} showed that the Hilbert scheme $H_{\mathbb P^3}^S$
of smooth connected curves in $\mathbb P^3$ is non-reduced.
$H_{\mathbb P^3}^S$ is the disjoint union of the open subscheme
$H_{d,g}^S$ consisting of curves of degree $d$ and genus $g$.
He considered a $56$-dimensional irreducible closed subset 
$W \subset H_{14,24}^S$
whose general member $C$ is contained in a smooth cubic surface.
He showed that the dimension of the tangent space of $H_{14,24}^S$
at $[C]$ is equal to $57$. Moreover, he proved that
$W$ is maximal as a subvariety of $(H_{14,24}^S)_{\red}$,
and hence $H_{14,24}^S$ is non-reduced.

We consider a generalization of Mumford's example.
Let $W$ be an irreducible closed subset of $H_{d,g}^S$
whose general member $C$ is contained in a smooth cubic surface.
Suppose that $W$ is maximal among all such subsets.
We ask the next question:
\begin{ques}\label{question}
 Is $W$ an irreducible component of $(H_{d,g}^S)_{\red}$?
 If so, is $H_{d,g}^S$ non-reduced at the generic point of $W$?
\end{ques}
\noindent
See \S 4 (\eqref{eqn:description of W} in particular)
for more explicit description of $W$.
First we observe that $\dim W=d+g+18$ when $d>9$,
while every irreducible component of $H_{d,g}^S$ is of dimension
at least $4d$.
Hence we consider the natural range
$\Omega:=\{(d,g)\in \mathbb Z^2|d > 9, g\ge 3d-18\}$ 
of pairs $(d,g)$, where the above question makes sense.
Secondly we consider the cohomology group
$H^1(\mathbb P^3,\mathcal I_C(3))$ for a general member $C$ of $W$.
The dimension $h^1(\mathcal I_C(3))$ as a vector space 
is the gap between $\dim W$ and
the dimension of the tangent space of $H_{d,g}^S$ at $[C]$.
This corresponds to the extra embedded
first order infinitesimal deformations of $C \subset \mathbb P^3$
other than the ones coming from $W$.
Thus if $h^1(\mathcal I_C(3))=0$, then
$W$ is an irreducible component of $H_{d,g}^S$
of $d+g+18$, and moreover, $H_{d,g}^S$ is non-singular at
the generic point of $W$.

In this paper, we concentrate on the case where $h^1(\mathcal I_C(3))=1$.
This is the first non-vanishing case, which includes Mumford's example.
In this case, there are only the two possibilities:
\begin{enumerate}
 \renewcommand{\labelenumi}{{\rm (\Alph{enumi})}}
 \item $H_{d,g}^S$ is non-reduced along $W$.
       Moreover, $W$ is an irreducible component of
       $(H_{d,g}^S)_{\red}$.
 \item There exists an irreducible component $V \supsetneqq W$ 
       of $H_{d,g}^S$ such that $\dim V=\dim W+1$ and 
       a general member is not contained in a cubic.
       Moreover, $H_{d,g}^S$ is generically smooth along $W$.
\end{enumerate}
We show that the case (B) does not occur.

\begin{thm}[Main Theorem]\label{thm:main}
 Let $(d,g) \in \Omega$ and let $W$ be an irreducible closed subset
 of $H_{d,g}^S$ whose general member $C$ is contained in a smooth 
 cubic surface.
 Suppose that $W$ is maximal among all such subsets.
 If $h^1(\mathcal I_C(3))=1$, then 
 $W$ is an irreducible component of $(H_{d,g}^S)_{\red}$
 of dimension $d+g+18$, and $H_{d,g}^S$ is non-reduced along $W$.
\end{thm}

For this kind of problem, two approaches are known.
One is to show that (B) leads to a contradiction, using e.g. liaison.
This was used by Mumford in \cite{Mumford}.
It has been also used to show that $H_{16,30}^S$ is non-reduced 
in \cite{Nasu}.
But it depends on case by case arguments.
Hence we cannot apply it for our general case 
that $h^1(\mathcal I_C(3))=1$.

In the proof of Theorem \ref{thm:main}, we use the other approach
described as follows. Let $C$ be a general member of $W$.
If $H_{d,g}^S$ is non-singular at $[C]$, then every
first order infinitesimal deformation $\varphi$
(i.e. a deformation over $\Spec k[t]/t^2$)
of $C \subset \mathbb P^3$ can be lifted to a deformation 
over $\Spec k[t]/(t^{n+1})$ for any integer $n\ge 2$.
We prove that there exists a first order infinitesimal deformation $\varphi$
of $C \subset \mathbb P^3$ that cannot be lifted to 
any deformation over $\Spec k[t]/(t^3)$
(cf. Proposition \ref{prop:core proposition}).
This implies that $H_{d,g}^S$ is singular along $W$, and hence we obtain (A).
This approach was first used by Curtin in \cite{Curtin},
who proved our result for the case of Mumford's example.
We generalize a calculation method used in his proof.
More precisely, we compute the obstruction map
$$
\begin{array}{ccc}
\Hom (\mathcal I_C,\mathcal O_C) 
& \rightarrow 
& \Ext^1(\mathcal I_C,\mathcal O_C) \\
\rotatebox{90}{$\in$} && \rotatebox{90}{$\in$} \\
\varphi & \mapsto & \varphi \cup \mathbf e \cup \varphi,
\end{array}
$$
where $\mathbf e \in \Ext^1(\mathcal O_C,\mathcal I_C)$ 
is the extension class of the basic exact sequence
\begin{equation}\label{ses:basic}
0 \longrightarrow \mathcal I_C \longrightarrow \mathcal O_{\mathbb P^3}
\longrightarrow \mathcal O_C \longrightarrow 0
\end{equation}
(cf. \S 2.1).
We use linear systems on the cubic surface $S$ containing $C$
for the computation.
Furthermore, we find an interesting
relation between the obstruction map and some geometry
arising from a conic pencil on the cubic $S$ (cf. \S 3.3).

Generalizations of Mumford's example were also studied
by Kleppe \cite{Kleppe85},\cite{Kleppe96} and Ellia \cite{Ellia}.
They gave a conjecture concerning non-reduced components of 
the Hilbert scheme $H_{\mathbb P^3}^S$ 
with some results which partially prove it
(see Conjecture \ref{conj:Kleppe-Ellia}).
Our theorem differently partially proves the conjecture.
See Remark \ref{rem:Kleppe-Ellia} for the relation between
their work and our theorem.
Constructions of non-reduced components of $H_{\mathbb P^3}^S$
by liaison or Rao module have been developed 
by Martin-Deschamps and Perrin \cite{MDP}, 
and by Fl{\o}ystad \cite{Floystad}.
See \cite{Floystad} for another generalization of Mumford's example.

{\bf Acknowledgements.}\quad
I should like to express my sincere gratitude
to my advisor, Professor Shigeru Mukai.
He read all the drafts of this paper very carefully,
pointed out a critical mistake in a draft,
and made many suggestions which greatly improved
the presentation and the proofs.
In particular, a discussion with him
led me to have the idea of using the Serre duality pairing
to improve a crucial part of the proof of
Proposition \ref{prop:core proposition}.
I am grateful to the referee for helpful comments.

\section*{Notation and Conventions}
We work in $\mathbb P^3$, the $3$-dimensional projective space over an
algebraically closed field $k$ of characteristic 0. 
Given a closed subscheme $V$ of $\mathbb P^3$,
we denote by $\mathcal I_V$ the ideal sheaf of $V$ in $\mathbb P^3$.
If $X \subseteq V$ is a closed subscheme, we indicate the ideal sheaf of $X$
in $V$ by $\mathcal I_{X/V}$. 
$\mathcal N_V \cong \sHom(\mathcal I_V,\mathcal O_V)$ 
and $\mathcal N_{X/V}$ 
denote the normal sheaf of $V$ in $\mathbb P^3$ and the normal sheaf 
of $X$ in $V$ respectively. 
Given $\mathcal O_{\mathbb P^3}$-modules $\mathcal F$ and $\mathcal G$, 
$h^i(\mathcal F)$, 
$\Hom(\mathcal F, \mathcal G)$ and
$\Ext^i(\mathcal F, \mathcal G)$
denote
$\dim H^i(\mathbb P^3, \mathcal F)$,
$\Hom_{\mathcal O_{\mathbb P^3}}(\mathcal F, \mathcal G)$ and
$\Ext^i_{\mathcal O_{\mathbb P^3}}(\mathcal F, \mathcal G)$
respectively. 
We denote the $p$-th {\v C}ech cohomology group of $\mathcal F$
with respect to an open covering $\mathfrak U$ by
$\check H^p(\mathfrak U,\mathcal F)$.
If $D$ is a Cartier divisor on a variety $X$, 
$\mathcal O_X(D)$ and $|D|$ respectively denote
the invertible sheaf and the complete linear system 
associated to $D$.
For a linear system $\Lambda$ on $X$, we denote the fixed
part of $\Lambda$ by $\Bs \Lambda$.
$\mathcal O_X(1)$ and $\mathbf h$ denote
the restriction of the tautological line bundle 
$\mathcal O_{\mathbb P^3}(1)$ to $X$ and the divisor class 
corresponding to $\mathcal O_X(1)$ respectively.
We denote by $\Rat (\mathcal L)$ the constant sheaf of global rational
sections of a line bundle $\mathcal L$ on $X$.
For a non-zero rational section $s$ of $\mathcal L$, 
we denote the divisor $(s)_0-(s)_{\infty}$ 
of zeros minus poles of $s$ by $\div(s)$.
%as $$\div(s):=\sum \ord_{\Gamma} (s) \Gamma $$
%in which $\ord_{\Gamma}(s)$ is the order of zeros or the minus of the
%order of poles of $s$ along a prime divisor $\Gamma \subset X$.
$\mathcal L(D)$ denotes the subsheaf of $\Rat(\mathcal L)$
which consists of rational sections $s$ of $\mathcal L$ such that
$\div(s)+D$ is effective.
We have $\mathcal L(D) \cong \mathcal L \otimes \mathcal O_X(D)$ 
by the usual multiplication map.

\section{Preliminaries}

\subsection{}
In this subsection, we recall some basic facts on the infinitesimal study 
of the Hilbert scheme of space curves. 
In what follows, we refer to  \cite[I.2]{Kollar} for the proofs,
where there is a very thorough discussion of general embedded deformations.

Let $C$ be a smooth connected curve in $\mathbb P^3$.
Then an ({\em embedded}) {\em $n$-th order} ({\em infinitesimal}) 
{\em deformation} 
of $C \subset \mathbb P^3$ is a closed subscheme 
$\mathcal C_n$ of $\mathbb P^3 \times \Spec k[t]/(t^{n+1})$
which is flat over $k[t]/(t^{n+1})$ and 
$\mathcal C_n \otimes_{k[t]/(t^{n+1})} k = C$.
The set of all first order deformations of $C \subset \mathbb P^3$
is the Zariski tangent space of $H_{\mathbb P^3}^S$ at the point $[C]$.
Let $\mathcal C_1$ be a first order deformation of 
$C \subset \mathbb P^3$. If there exists no second order deformation 
$\mathcal C_2$ of $C \subset \mathbb P^3$
such that $\mathcal C_2 \otimes_{k[t]/(t^3)} k[t]/(t^2) = \mathcal C_1$,
we say $\mathcal C_1$ is {\em obstructed at the second order}.
The set of all first order deformations of $C \subset \mathbb P^3$
is parametrized by $\Hom(\mathcal I_C,\mathcal O_C)$.
So we abusively identify them from now.
The basic exact sequence \eqref{ses:basic} induces the isomorphism
$$
\delta:\Hom(\mathcal I_C,\mathcal O_C) 
\overset{\sim}{\longrightarrow} \Ext^1(\mathcal I_C,\mathcal I_C).
$$ 
Let $\varphi \in \Hom(\mathcal I_C,\mathcal O_C)$ 
be a first order deformation of $C \subset \mathbb P^3$.
Then $\varphi$ is obstructed at the second order
if and only if the cup product 
$o(\varphi):=\delta(\varphi)\cup \varphi$ by
\begin{equation}\label{map:cup product1}
\cup_1: \Ext^1(\mathcal I_C, \mathcal I_C)
\times \Hom(\mathcal I_C, \mathcal O_C) 
\overset \cup\longrightarrow \Ext^1(\mathcal I_C,\mathcal O_C)
\end{equation}
is non-zero. $o(\varphi)$ is called the {\em obstruction} to extend
$\varphi$ to second order deformations.
Since $C$ and $\mathbb P^3$ are both non-singular,
$C$ is a local complete intersection in $\mathbb P^3$.
Therefore the obstruction $o(\varphi)$ 
is contained in $H^1(\mathcal N_C)$, which is regarded as a subspace of 
$\Ext^1(\mathcal I_C,\mathcal O_C)$ by
the exact sequence
$$
0 \longrightarrow H^1(\sHom(\mathcal I_C, \mathcal O_C))
\longrightarrow \Ext^1(\mathcal I_C, \mathcal O_C)
\longrightarrow H^0(\sExt^1(\mathcal I_C, \mathcal O_C))
\longrightarrow 0
$$
obtained from local-global spectral sequence for $\Ext$.

From now on, we treat the case where
$C$ is contained in a smooth cubic surface $S$.
The natural sheaf inclusion
$\mathcal O_{\mathbb P^3}(-3) \cong 
\mathcal I_S \overset{\iota}{\hookrightarrow} \mathcal I_C$ induces 
the homomorphisms
\begin{align*}
\psi &:\Hom(\mathcal I_C,\mathcal O_C) \longrightarrow 
\Hom(\mathcal I_S,\mathcal O_C) \cong H^0(\mathcal O_C(3)), \\
\psi'&:\Ext^1(\mathcal I_C,\mathcal I_C) \longrightarrow 
\Ext^1(\mathcal I_S,\mathcal I_C) \cong H^1(\mathcal I_C(3)),
\qquad \mbox{and} \\
\psi''&:\Ext^1(\mathcal I_C,\mathcal O_C) \longrightarrow 
\Ext^1(\mathcal I_S,\mathcal O_C) \cong H^1(\mathcal O_C(3)).
\end{align*}
We denote by $\pi$ the composite
\begin{equation}\label{map:obstruction}
\psi''\circ o = \psi'' \circ \cup_1 \circ (\delta \times \id):
\Hom(\mathcal I_C,\mathcal O_C) \rightarrow H^1(\mathcal O_C(3)). 
\end{equation}
Then the following is obvious.

\begin{prop}[\cite{Floystad} Corollary 1.3]\label{prop:criterion}
 Let $\varphi$ be an embedded first order infinitesimal 
 deformation of a curve $C \subset \mathbb P^3$ 
 on a smooth cubic surface $S$. 
 If $\pi(\varphi)$ is non-zero in $H^1(\mathcal O_C(3))$,
 then $\varphi$ is obstructed at the second order. 
 \qed
\end{prop}

Let us give another expression of $\pi$.
A natural cup product map
\begin{equation}\label{map:cup product2}
\cup_2: H^1(\mathcal I_C(3))
\times \Hom(\mathcal I_C, \mathcal O_C) 
\overset \cup\longrightarrow H^1(\mathcal O_C(3))
\end{equation}
satisfies the commutative diagram
$$
\begin{array}{cccc}
   \Ext^1(\mathcal I_C, \mathcal I_C) \times \Hom(\mathcal I_C,\mathcal O_C) &
   \mapright{\cup_1} & \Ext^1(\mathcal I_C,\mathcal O_C) \\
  \mapdown{(\psi',\id)} && \mapdown{\psi''} \\
   H^1(\mathcal I_C(3)) \times \Hom(\mathcal I_C,\mathcal O_C) &
   \mapright{\cup_2} & 
   H^1(\mathcal O_C(3)). \\
\end{array}
$$
Moreover, $\psi$ and $\psi'$ naturally satisfy a commutative diagram
$$
\begin{array}{ccc}
 \Hom(\mathcal I_C,\mathcal O_C) &
  \mapright{\delta} & \Ext^1(\mathcal I_C,\mathcal I_C) \\
 \mapdown{\psi} && \mapdown{\psi'} \\
 H^0(\mathcal O_C(3)) &
 \mapright{\bar \delta} & H^1(\mathcal I_C(3)), \\
\end{array}
$$
where $\bar \delta$ is the coboundary map of 
\begin{equation}\label{les:cubic normality}
0 \longrightarrow H^0(\mathcal I_C(3)) \longrightarrow
H^0(\mathcal O_{\mathbb P^3}(3)) \longrightarrow
H^0(\mathcal O_C(3)) \overset {\bar \delta}\longrightarrow
H^1(\mathcal I_C(3)) \longrightarrow
0
\end{equation}
induced from \eqref{ses:basic}$\otimes\mathcal O_{\mathbb P^3}(3)$.
Hence we have another expression of $\pi$ as
\begin{equation}\label{eqn:obst map}
\pi
=\psi'' \circ \cup_1 \circ (\delta \times \id)
=\cup_2 \circ ((\psi'\circ \delta)\times \id)
=\cup_2 \circ ((\bar\delta \circ \psi)\times \id). 
\end{equation}
By definition, $\psi$ maps an element $\varphi$ of 
$\Hom(\mathcal I_C,\mathcal O_C)$ to $u=\varphi_3(f)$,
where $f$ is the cubic polynomial which defines the
isomorphism $\mathcal O_{\mathbb P^3}(-3) \cong \mathcal I_S$,
and $\varphi_3$ is the homomorphism 
$H^0(\mathcal I_C(3)) \rightarrow H^0(\mathcal O_C(3))$ 
induced from $\varphi$. Moreover, $\psi$ is surjective.

\subsection{}
In this subsection, we recall some basic facts on 
linear systems on a smooth cubic surface.
Let $\mathcal L$ be an invertible sheaf on a smooth cubic surface $S$.
We may consider $S$ to be a
$\mathbb P^2$ blown up at 6 points in a
general position and embedded
by anti-canonical linear system $|-K_S|$ in $\mathbb P^3$. 
The classes of the pull back $\mathbf l$ of a line in $\mathbb P^2$ and six 
exceptional curves $\mathbf e_i$ ($1 \le i \le 6$) 
form a $\mathbb Z$-free basis of the Picard group $\Pic S$ of $S$.
Thus there is an isomorphism $\Pic S \cong {\mathbb Z}^{\oplus 7}$
sending the class
$\mathcal L = a \mathbf l - \sum_{i=1}^6 b_i \mathbf e_i$
to a 7-tuple $(a;b_1,\ldots,b_6)$ of integers.
We denote the class $3\mathbf l -\sum_{i=1}^6 \mathbf e_i$ 
of hyperplane sections by $\mathbf h$.
Recall that the Weyl group $W(\mathbb E_6)$ acts on $\Pic S$.
By virtue of this action, we can choose a suitable blow-up
$S \rightarrow \mathbb P^2$ for $\mathcal L$ such that
\begin{equation}\label{eqn:W-standard}
 b_1 \ge b_2 \ge \ldots \ge b_6 \quad {\rm and} \quad a\ge b_1+b_2+b_3
\end{equation}
holds. 
When \eqref{eqn:W-standard} holds, we say the basis 
$\left\{ \mathbf l,\mathbf e_1, \ldots, \mathbf e_6 \right\}$
is {\em $\mathbb E$-standard} for $\mathcal L$.
The 7-tuple $(a;b_1,\ldots,b_6)$ 
is uniquely determined for each invertible sheaf $\mathcal L$ on $S$.
We call it the {\em $\mathbb E$-multidegree} of $\mathcal L$.
For a divisor $D$ on $S$, we define the $\mathbb E$-multidegree of $D$ 
as that of the associated invertible sheaf $\mathcal O_S(D)$.

$\mathbb E$-standard basis is useful for analyzing the linear system
$|D|$ associated to a divisor $D$ on a smooth cubic surface.
\begin{lem}\label{lem:zariski decomp}
 Let $D$ be a divisor of $\mathbb E$-multidegree $(a;b_1,\ldots,b_6)$
 on a smooth cubic surface $S$.
 \begin{enumerate}
  \renewcommand{\labelenumi}{{\rm (\roman{enumi})}}
  \item The following are equivalent:
	\begin{enumerate}
	 \item $D\ge 0$ and $|D|$ is $($base point$)$ free;
	 \item $D$ is nef $($i.e. $D\cdot C\ge 0$ for any curve $C$ on $S)$;
	 \item $b_6 \ge 0$.
	\end{enumerate}
  \item If $b_6 \ge 0$, then $D^2 \ge 0$.
	The equality holds if and only if $a=b_1$.
  \item If $|D|\ne \emptyset$,
	then the fixed part of $|D|$ is
	$$F=\sum\limits^6_{\substack{i=1 \\ b_i < 0}} 
	(-b_i) \mathbf e_i$$ 
	for an $\mathbb E$-standard basis for $\mathcal O_S(D)$.
 \end{enumerate}
\end{lem}
Here we abusively identify the class $\mathbf e_i$
with the unique effective divisor in the class.
We refer to Geramita \cite{Geramita} for the proof.

\vskip 3mm

When $C$ is a smooth connected curve on a smooth cubic surface,
the $\mathbb E$-multidegree $(a;b_1,\ldots,b_6)$ of $C$ 
satisfies $b_6 \ge 0$ if $C$ is not a line,
and $a > b_1$ if $C$ is not a conic.

Let $F$ be a ``multiple line'' or a ``multiple conic'' 
on a smooth cubic surface.
We compute $h^i$ ($i=0,1$) of the structure sheaf of $F$:
\begin{lem}\label{lem:multiple curve}
 Let $m>0$ and let $mE$ $($resp. $mD)$ be a member of 
 the linear system $|m\mathbf e_1|$ 
 $($resp. $|m(\mathbf l -\mathbf e_1)|)$
 on a smooth cubic surface $S$. Then we have
$$
 \begin{array}{lcc}
  \dim H^0(\mathcal O_{mE}) = \dfrac{m(m+1)}2, 
   & 
   & H^1(\mathcal O_{mE})=0, \\
  \dim H^0(\mathcal O_{mD}) = m, 
   & \mbox{and} 
   & H^1(\mathcal O_{mD})=0. \\
 \end{array}
 $$
\end{lem}
\noindent
{\sc Proof. \quad}
We prove the assertion for a multiple line $mE$ by induction on 
$m \in \mathbb N$. It is clear for $m=1$.
There exists an exact sequence
$$
0 \longrightarrow \ker q
\longrightarrow \mathcal O_{mE}
\overset{q}{\longrightarrow} \mathcal O_{(m-1)E}
\longrightarrow 0.
$$
Since the sheaf $\bar {\mathcal I}_{mE}$
of ideal defining $mE$ in $S$ is isomorphic to
$\mathcal O_S(-mE)$, we have isomorphisms
$$
\ker q 
\cong \bar {\mathcal I}_{mE} / \bar {\mathcal I}_{(m-1)E}
\cong \mathcal O_S(-(m-1)E)\big{\vert}_E
\cong \mathcal O_{\mathbb P^1}(m-1).
$$
Therefore, by the inductive assumption, we get
$$
h^0(\mathcal O_{mE})
=h^0(\ker q)+h^0(\mathcal O_{(m-1)E})
=m+m(m-1)/2=m(m+1)/2
$$ and
$H^1(\mathcal O_{mE})=0$.
The proof for a multiple conic $mD$ is similar (use $D^2=0$).
\qed

\vskip 3mm

We next characterize the freeness of $|D|$ for a divisor $D$
by the vanishing of $H^1(S,-D)$.
Let $D$ be a non-zero effective divisor on a smooth cubic surface $S$.
Then $|D|$ has the unique decomposition
$$|D|=|D'|+F,$$ 
where $F:=\Bs |D|$ and $|D'|$ is free by Lemma \ref{lem:zariski
decomp} (i). 
(When $|D|$ is free, $D=D'$ and $F=0$.)
\begin{lem}\label{lem:vanishing thm}
 Let $D$, $D'$, and $F$ be as above. Then
 \begin{enumerate}
  \item We have
	$$
	h^1(S,-D)=h^0(\mathcal O_{D'})+h^0(\mathcal O_F)-1.
	$$
	If $(D')^2>0$, then
	a general member of $|D'|$ is 
	a smooth connected curve and hence $h^0(\mathcal O_{D'})=1$.
	If $(D')^2=0$, then a general member of $|D'|$ is 
	a disjoint union of $m$ conics 
	for some $m \in \mathbb Z_{\ge 0}$ and hence
	$h^0(\mathcal O_{D'})=m$.
  \item Suppose that $D^2>0$. Then we have
	$h^1(S,-D)=h^0(\mathcal O_F)$.
	In particular, $|D|$ is free 
	if and only if $H^1(S,-D)=0$.\footnote{
	The only-if part is a particular consequence of
	Kawamata-Viehweg's vanishing theorem that
	$H^i(X,K_X+D)=0$ for a nef (i.e. $D\cdot C\ge 0$ for any curve $C$)
	and big (i.e. $D^2>0$) Cartier divisor $D$ on a smooth
	surface $X$ and for $i>0$.
	In what follows, we say ``$D$ is nef and big'' to mean
	that $|D|$ is free and $D^2>0$.}
  \item If $|D|$ is free, then $H^i(S,D)=0$ for $i=1,2$.
 \end{enumerate}
\end{lem}
\noindent
{\sc Proof. \quad}
(1) Let $D$, $D'$, and $F$ be as above. Since $D$ is effective, we have
$h^1(S,-D)=h^0(\mathcal O_D)-1$
by an exact sequence 
$0 \rightarrow \mathcal O_S(-D) 
\rightarrow \mathcal O_S
\rightarrow \mathcal O_D
\rightarrow 0$.
By Lemma \ref{lem:zariski decomp} (iii), 
$D'$ and $F$ have disjoint supports.
Therefore, we have
$\mathcal O_D \cong \mathcal O_{D'} \oplus \mathcal O_F$
and
$h^1(S,-D)=h^0(\mathcal O_{D'})+h^0(\mathcal O_F)-1$.
When $(D')^2>0$, 
$D'$ is ample or a pull-back of an ample divisor on
a $\mathbb P^2$ blown-up at less than $6$ points.
Therefore, a general member of $|D'|$ is a smooth connected curve
by Bertini's theorem.
When $(D')^2=0$, $D'$ is linearly equivalent to 
$m(\mathbf l -\mathbf e_1)$ for some $m\ge 0$ by Lemma
\ref{lem:zariski decomp} (ii), which is the class of
$m$ conics. Therefore,
the case is also a consequence of Bertini's theorem
together with Lemma \ref{lem:multiple curve}.

(2) Let $F$ be the fixed part of $|D|$.
Then $F$ is a disjoint sum of (multiple) lines or zero.
Thus we have $F^2 \le 0$.
Since $D'(\sim D-F)$ and $F$ are disjoint, we get $D\cdot F=F^2$.
Therefore $D^2>0$ implies $(D')^2>0$
by $(D')^2=(D-F)^2=D^2-2D\cdot F+F^2=D^2-F^2$.
If $F \ne 0$, then we get
$h^1(S,-D)=h^0(\mathcal O_F)\ne 0$ by (1).
If $F =0$ (i.e. $|D|$ is free), then we get
$h^1(S,-D)=h^0(\mathcal O_D)-1=0$.

(3) Let $\mathbf h$ be the class of hyperplane sections of $S$.
Since $|D|$ is free, $D+\mathbf h$ is very ample.
By the Serre duality and the Kodaira vanishing theorem, we have
$H^i(S,D)\cong H^{2-i}(S,-(D+\mathbf h))=0$ for $i=1,2$.
\qed

\vskip 3mm

We use Lemma \ref{lem:vanishing thm} to compute 
$h^1(\mathcal I_C(n))$ ($n \in \mathbb Z$)
for a curve $C$ on a smooth cubic surface $S$.
Let $(a;b_1,\ldots,b_6)$ be the $\mathbb E$-multidegree of $C$.
Given $n \in \mathbb Z_{\ge 0}$,
we consider the linear system $\Lambda_n:=|C-n\mathbf h|$ on $S$,
where $\mathbf h=(3;1,\ldots,1)$ 
is the class of hyperplane sections.
Suppose that $\Lambda_n \ne \emptyset$.
Then by Lemma \ref{lem:zariski decomp} (iii),
the fixed part $F$ of $\Lambda_n$
is a disjoint sum of (multiple) lines as follows:
\begin{equation}\label{form:fixed parts}
F=\sum^6_{\substack{i=1 \\ b_i < n}}
F_i, \quad F_i := (n-b_i)E_i
\end{equation}
for an $\mathbb E$-standard basis for $C$.
Here each $E_i$ ($1 \le i \le 6$) denotes the line
corresponding to the class $\mathbf e_i$ of exceptional curve.
Since all $F_i$'s are disjoint, we have
$\mathcal O_F \cong \bigoplus_{b_i <n} \mathcal O_{F_i}$.
By Lemma \ref{lem:multiple curve}, we get 
$h^0(\mathcal O_{F_i})=(n+1-b_i)(n-b_i)/2$ for every $i$.
The exact sequence 
$
0 \rightarrow \mathcal I_S(n) 
\rightarrow \mathcal I_C(n)
\rightarrow \mathcal I_{C/S}(n)
\rightarrow 0
$
induces an isomorphism
\begin{equation}\label{isom:rao component}
H^1(\mathcal I_C(n)) 
\cong H^1(\mathcal I_{C/S}(n))
\cong H^1(S,-(C-n\mathbf h)). 
\end{equation}
Thus we have the next corollary by applying Lemma
\ref{lem:vanishing thm} (2) to $D=C-n\mathbf h$.

\begin{cor}\label{cor:vanishing thm2}
 Let $C$ be a smooth connected curve of $\mathbb E$-multidegree
 $(a,b_1,\ldots,b_6)$ on a smooth cubic surface $S$.
 Assume that $\Lambda_n:=|C-n\mathbf h|\ne \emptyset$ 
 and $(C-n\mathbf h)^2>0$ for $n \in \mathbb Z_{\ge 0}$.
 Then we have
 $$
 h^1(\mathcal I_C(n)) 
 = h^0(\mathcal O_F)  
 =\sum_{\substack{i=1 \\ b_i < n}}^6 \frac{(n+1-b_i)(n-b_i)}2,
 $$
 where $F=\Bs \Lambda_n$.
 In particular, $\Lambda_n$ is free if and only 
 if $H^1(\mathcal I_C(n))=0$.
\end{cor}

\subsection{}
In this subsection, we define some restriction maps.
In what follows, 
when $X$ is a subscheme of $\mathbb P^3$ and $F$ is a polynomial 
of degree $d$, we sometimes use the same symbol $F$ 
to denote the element $F\big{\vert}_X$ of $H^0(\mathcal O_X(d))$
if there is no confusion.

Let $S$ be a smooth cubic surface and let $E$ be a line on $S$.
Let $x,y$ be two linear forms on $\mathbb P^3$ defining $E$. 
Then the cubic polynomial $f$ defining $S$ is
\begin{equation}\label{eqn:cubic}
 f=Ax+By
\end{equation}
for two quadratic polynomials $A,B$ on $\mathbb P^3$.
By definition, $x,y$ form a basis of $H^0(\mathcal O_S(1)(-E))$.
The corresponding linear system $\Lambda=|\mathbf h - E|$ 
defines the projection $p:S \rightarrow \mathbb P^1$ from $E$.
By this map, $S$ has a conic bundle structure.
Let $x',y'$ be the sections of $p^{\ast}\mathcal O_{\mathbb P^1}(1)$
corresponding to $x,y$. Then $S$ is covered by two open subsets
$D(x')$ and $D(y')$ of $S$. Let $s$ be a rational section of 
$\mathcal O_S(1)$ defined by
\begin{equation}\label{eqn:def of s}
 s =
  \begin{cases}
   \dfrac {-B}x & \mbox{on $D(x')$}, \vspace{1mm} \\
   \dfrac Ay & \mbox{on $D(y')$}.
  \end{cases}
\end{equation}
Then by construction, $s$ is a global section of $\mathcal O_S(1)(E)$.
Moreover, by the correspondence
$$
\begin{array}{ccc}
m: \mathcal O_E &\overset{\sim}{\longrightarrow}
 & \mathcal O_S(1)(E)\big{\vert}_E \\
\phantom{m:} \rotatebox{90}{$\in$} && \rotatebox{90}{$\in$} \\
\phantom{m:} \mu & \longleftrightarrow & \mu (s\big{\vert}_E),
\end{array}
$$
we get a trivialization of the line bundle
$\mathcal O_S(1)(E)\big{\vert}_E \cong \mathcal O_{\mathbb P^1}$.
Applying $\otimes \mathcal O_E(n-1)$ to $m^{-1}$,
we have a natural isomorphism
$m_n:\mathcal O_S(n)(E)\big{\vert}_E 
\overset{\sim}{\rightarrow} \mathcal O_E(n-1)$.
We define a homomorphism $r_E$ of
$\mathcal O_S$-modules by the composite
$$
r_E:\mathcal O_S(n)(E) 
\overset{\res}{\longrightarrow} \mathcal O_S(n)(E)\big{\vert}_E
\overset{m_n}{\longrightarrow} \mathcal O_E(n-1),
$$
where $\res$ is the restriction map.
Then we have an exact sequence 
\begin{equation}\label{ses:restriction map}
 0 
\longrightarrow \mathcal O_S(n)
\longrightarrow \mathcal O_S(n)(E)
\overset{r_E}{\longrightarrow} \mathcal O_E(n-1)
\longrightarrow
0.
\end{equation}

We explicitly describe the restriction map $H^0(r_E)$
for any positive integer $n$.
Let $v$ be an element of $H^0(\mathcal O_S(n)(E))$.
Then the multiplication map
$$
H^0(\mathcal O_S(1)(-E)) 
\otimes 
H^0(\mathcal O_S(n)(E))
\longrightarrow
H^0(\mathcal O_S(n+1))
$$
gives two elements $x v,y v$ in $H^0(\mathcal O_S(n+1))$.
Since $H^0(\mathcal O_{\mathbb P^3}(n+1)) 
\rightarrow H^0(\mathcal O_S(n+1))$ is surjective,
there exist two polynomials $\eta_1,\eta_2$ of degree $n+1$ such that
their restrictions to $S$ are $x v$, $y v$.
Hence we have an equality
\begin{equation}\label{eqn:on S}
v =\dfrac{\eta_1}x =\dfrac{\eta_2}y
\quad \mbox{in} \quad \Rat(\mathcal O_S(n)). 
\end{equation}
Since $v$ is globally defined,
there exists a polynomial $\xi$ of degree $n-1$ such that
\begin{equation}\label{eqn:relation}
 x\eta_2 - y \eta_1 = \xi f. 
\end{equation}
Here we see that $\xi\big{\vert}_E$ does not depend on 
the choice of $\eta_1,\eta_2$.
Here and later, for a polynomial $F$, we denote 
$F \pmod{\langle x, y \rangle}$ by $F \big{\vert}_E$.
We show that $\xi\big{\vert}_E$ agrees with $r_E(v)$.

\begin{claim}\label{claim:explicit description of r_E}
 $r_E(v)=\xi \big{\vert}_E$.
\end{claim}
\noindent
{\sc Proof. \quad}
Since $f=Ax+By$, we have
$x(\eta_2-A\xi)=y(\eta_1+B\xi)$ by \eqref{eqn:relation}.
Since $x$ and $y$ are coprime,
there exists a polynomial $\eta'$ of degree $n$
such that $\eta_1=-B\xi+x\eta'$ and $\eta_2=A\xi+y\eta'$.
Therefore, we obtain $v=\xi s+\eta'$ from
\eqref{eqn:on S} and \eqref{eqn:def of s}.
We see $r_E(\eta')=0$ because $\eta'$ is a polynomial.
Hence we get $r_E(v)=r_E(\xi s+\eta')=\xi\big{\vert}_E$
from $(\xi s)\big{\vert}_E=(\xi\big{\vert}_E) (s\big{\vert}_E)$
and the trivialization $m$.
\qed

Thus we get the description of $H^0(r_E)$.

\begin{rem}\label{rem:conic pencil}
Let $\Lambda$ be the linear system $|\mathbf h -E|$ 
corresponding to $\mathcal O_S(1)(-E)$.
Then the restriction $\Lambda\big{\vert}_E$ is a subpencil of
$|\mathcal O_E(2)| \cong |\mathcal O_{\mathbb P^1}(2)|$ 
since $(\mathbf h-E)\cdot E=2$.
Writing the cubic equation $f$ in the form $f=Ax+By$ 
is also useful to describe the restriction map in this case.
We see that planes $H$ through $E$ are parametrized 
by $\mathbb P^1_{(t_0,t_1)}$ and $H=H_{(t_0,t_1)}$
defined by $t_0 x +t_1 y=0$.
A member of $\Lambda$ is a conic
defined by $t_0 x +t_1 y=t_0 (-B)+t_1 A=0$.
Hence a member of $\Lambda\big{\vert}_E$ is a divisor of 
degree two on $E$, which is defined by 
$t_0 (-B)\big{\vert}_E+t_1 A\big{\vert}_E=0$.

By a similar argument, we have a natural isomorphism
$\mathcal O_S(1)(-E)\big{\vert}_E \cong \mathcal O_E(2)$.
The composition of the restriction map
$\mathcal O_S(1)(-E) \overset{\res}{\rightarrow} O_S(1)(-E)\big{\vert}_E$
and the isomorphism induces
$$
r_E: H^0(\mathcal O_S(1)(-E)) \longrightarrow H^0(\mathcal O_E(2)),
$$ 
which sends
$t_0 x +t_1 y$ to $t_0 (-B)\big{\vert}_E + t_1 A\big{\vert}_E$.
We can see the one-to-one correspondence between
$|\im r_E|$ and $\Lambda \big{\vert}_E$
by taking the divisor of zeros.
\end{rem}

\section{Obstructed deformation of space curves}

We devote the whole section to the proof of the next proposition.

\begin{prop}[Core Proposition]\label{prop:core proposition}
 Let $S$ be a smooth cubic surface, 
 let $\mathbf h$ be the class of hyperplane sections,
 and let $\mathbf D$ be a divisor class of $S$ satisfying
 \begin{enumerate}
  \renewcommand{\labelenumi}{{\rm (\roman{enumi})}}
  \item The fixed part of the linear system $|\mathbf D-3\mathbf h|$ on $S$ is
        exactly a line $E$,
  \item $|\mathbf D-4\mathbf h|\ne \emptyset$.
 \end{enumerate}
 Then any general member $C$ of $|\mathbf D|$ has some
 embedded first order infinitesimal deformation
 which is obstructed at the second order.
\end{prop}

First we observe $|\mathbf D|\ne \emptyset$ by (ii).
Moreover, since both $|\mathbf D -3\mathbf h -E|$ and $|3\mathbf h +E|$ 
are free by assumption and Lemma \ref{lem:zariski decomp} (i),
a general member $C$ of $|\mathbf D|$ is a smooth connected curve
by Bertini's theorem.
Let $S$, $\mathbf h$, $\mathbf D$, $E$, and $C$ be as in the statement.
Let $x$, $y$, $A$, $B$ and $f$ be as in \S 2.3.
We fix these notation throughout the proof.
Now we start the proof.

\begin{claim}\label{claim:finite subscheme $Z$}
Let $Z:=C\cap E$, then $Z$ is of length two.
\end{claim}
\noindent
{\sc Proof. \quad}
Let $(a;b_1,\ldots,b_6)$ be the multidegree of $C$ on $S$
and let $\left\{\mathbf l,\mathbf e_1,\ldots,\mathbf e_6 \right\}$
be an $\mathbb E$-standard basis of $\Pic S$ for $C$.
Then by Lemma \ref{lem:zariski decomp} (iii),
the fixed part $\Bs|C-3\mathbf h|$ is a sum 
$\sum (3-b_i)\mathbf e_i$ over all $b_i < 3$.
On the other hand, we have $\Bs|C-3\mathbf h|=E$ by assumption.
Hence we have $E=\mathbf e_6$ and $b_6=2$.
This implies $C\cdot E=b_6=2$.
\qed

\vskip 3mm

\begin{lem}\label{lem:general $Z$}
 Let $\Lambda$ be the conic pencil $|\mathbf h -E|$ on $S$
 and let $\Lambda \big{\vert}_E$ be its restriction to $E$. 
 $($We refer to Remark $\ref{rem:conic pencil}$.$)$
 Then, $Z$ is not a member of $\Lambda \big{\vert}_E$.
\end{lem}
\noindent
{\sc Proof. \quad}
There exists an exact sequence 
$$
0 \longrightarrow \mathcal O_S(\mathbf D-E)
\longrightarrow \mathcal O_S(\mathbf D)
\longrightarrow \mathcal O_S(\mathbf D)\big{\vert}_E 
\longrightarrow 0.
$$
Then Lemma \ref{lem:vanishing thm} (3) shows
$H^1(S,\mathbf D-E)=0$ because $|\mathbf D-E|$ is free.
Hence the restriction map
$H^0(\mathcal O_S(\mathbf D)) \rightarrow H^0(\mathcal O_S(\mathbf D)\big{\vert}_E)$
is surjective.
We know $\dim |\mathcal O_S(\mathbf D)\big{\vert}_E|=2$
by Claim \ref{claim:finite subscheme $Z$}, while
$\Lambda \big{\vert}_E$ is a pencil. Thus we have
$\Lambda \big{\vert}_E \subsetneqq 
|\mathcal O_S(\mathbf D)\big{\vert}_E|$.
Therefore, any general member $C$ of $|\mathbf D|$ 
meets $E$ at $Z \notin \Lambda \big{\vert}_E$.
\qed

\begin{claim}\label{claim:nef and big}
   $C-3\mathbf h - E$ is nef and big.
\end{claim}
\noindent
{\sc Proof. \quad}
Put $D:=C-3\mathbf h - E$.
Since $D$ is cleary nef by assumption,
it suffices to show $D^2>0$.
Put $D_1:=C-4\mathbf h$ and $D_2:=\mathbf h-E$.
Then $D=D_1+D_2$. We obtain
$$
D^2 \ge D\cdot D_2 = (D_1+D_2)\cdot D_2=D_1\cdot D_2,
$$
because $D$ is nef, $D_1$ is effective, and $(D_2)^2=0$.
Since $D_1\cdot D_2
=(C-4\mathbf h)\cdot (\mathbf h-E)
\ge (C-4\mathbf h)\cdot (-E)=-2+4=2$, we have $D^2>0$.
\qed

\vskip 3mm

Since $\Bs|C-3\mathbf h|=E$ and $(C-3\mathbf h-E)^2>0$,
we have $h^1(S,-(C-3\mathbf h))=1$ by Lemma \ref{lem:vanishing thm} (1).
Hence we get $h^1(\mathcal I_C(3))=1$ by \eqref{isom:rao component}.
Thus there exists an element $u$ of $H^0(\mathcal O_C(3))$
which is not (the image of) a cubic polynomial, and 
an element $\varphi$ of $\Hom(\mathcal I_C,\mathcal O_C)$
such that $\varphi(f)=u$ (cf. the last paragraph of \S 2.1).
Let $\pi$ be the map defined by \eqref{map:obstruction}.
Then we have
$\pi (\varphi)
= (\bar \delta (\psi (\varphi))) \cup \varphi
= (\bar \delta (\varphi(f))) \cup \varphi
= \bar \delta (u) \cup \varphi$
by the alternative expression \eqref{eqn:obst map} of $\pi$.
Thus it suffices to show the following:
the cup product $\bar \delta(u) \cup \varphi$ by
$$
\cup_2:H^1(\mathcal I_C(3)) \times \Hom(\mathcal I_C,\mathcal O_C)
\overset{\cup}{\longrightarrow} H^1(\mathcal O_C(3))
$$
is non-zero in $H^1(\mathcal O_C(3))$.
(See \S 2.1 for $\bar \delta,\cup_2$ etc.)
If it is proved, then by Proposition \ref{prop:criterion}, 
$\varphi$ is obstructed at the second order.
Our procedure for this is as follows:
we relate the above cup
product map to familiar Serre duality pairing 
via several cup product maps, 
and eventually obtain the non-zero of the original product
from the perfect pairing.
First of all, since 
$\Hom(\mathcal I_C,\mathcal O_C)
\cong H^0(\sHom_{\mathbb P^3} (\mathcal I_C, \mathcal O_C))
\cong H^0(\mathcal N_C)$
is a cohomology group on $C$, the above $\cup_2$ is
compatible with the cup product map
$$
\cup_3:H^1({\mathcal N_C}^{\vee}(3)) \times H^0(\mathcal N_C)
\overset{\cup}{\longrightarrow} H^1(\mathcal O_C(3))
$$
via natural maps. Here ${\mathcal N_C}^{\vee}$ is the conormal bundle
$\mathcal I_C/{\mathcal I_C}^2$ of $C$. Moreover,
since $Z$ is an effective divisor on $C$,
by tensoring $\mathcal O_C(2Z)$ with the first and the last sheaves
of $\cup_3$, we get another cup product map
$$
\cup_4:H^1(({\mathcal N_C}^{\vee}(3)(2Z)) \times H^0(\mathcal N_C)
\overset{\cup}{\longrightarrow} H^1(\mathcal O_C(3)(2Z)),
$$
which is also compatible with the previous ones $\cup_i$ ($i=1,2,3$)
via natural maps.
\vskip 3mm

\subsection{}
In this subsection, we compute the obstruction. Let $u$ be as above.
By the exact sequence \eqref{ses:restriction map} as $n=3$,
we have a commutative diagram of exact sequences
\begin{equation}\label{diag:residue diagram}
\begin{array}{ccccccccc}
&& 0 && 0 && 0 && \\
&& \downarrow && \downarrow && \downarrow && \\
0 & \rightarrow & \mathcal O_S(3)(-C)
& \mapright{} & \mathcal O_S(3)(E-C)
& \mapright{r_E} & \mathcal O_E(2)(-Z)
& \rightarrow & 0 \\
&& \mapdown{} && \mapdown{} && \mapdown{} && \\
0 & \rightarrow & \mathcal O_S(3)
& \mapright{} & \mathcal O_S(3)(E)
& \mapright{r_E} & \mathcal O_E(2)
& \rightarrow & 0 \\
&& \mapdown{\res} && \mapdown{} && \mapdown{} && \\
0 & \rightarrow & \mathcal O_C(3)
& \mapright{} & \mathcal O_C(3)(Z)
& \mapright{} & \mathcal O_Z
& \rightarrow & 0 \\
&& \downarrow && \downarrow && \downarrow && \\
&& 0 && 0 && 0. && 
\end{array}
\end{equation}
Since $C-3\mathbf h-E$ is nef and big 
by Claim \ref{claim:nef and big}, we have
$$H^i(\mathcal O_S(3)(E-C))=0  \qquad (i=0,1)$$
by Lemma \ref{lem:vanishing thm} (2).
Hence the diagram induces an isomorphism
$
H^0(\mathcal O_S(3)(E)) 
\overset{\sim}{\rightarrow} H^0(\mathcal O_C(3)(Z)).
$
Thus there exists an element 
$\hat u$ of $H^0(\mathcal O_S(3)(E))$ 
such that $\hat u \big{\vert}_C=u$.
In particular, as we saw in \S 2.3
(cf. \eqref{eqn:on S} and \eqref{eqn:relation}),
there exist a quadratic polynomial $\xi$
and two quartic polynomials $\eta_1, \eta_2$ such that
\begin{equation}\label{presentation of $u$}
u=\dfrac{\eta_1}x=\dfrac {\eta_2}y \quad \mbox{in} \quad \Rat(\mathcal O_C(3))
\quad \mbox{and} \quad x\eta_2 -y\eta_1 = \xi f \quad \mbox{as a polynomial}.
\end{equation}
Moreover, by the snake lemma, we have
$$
H^1(\mathcal O_S(3)(-C)) 
=\Coker H^0(\res) \overset{\sim}{\longrightarrow} H^0(\mathcal O_E(2)(-Z)).
$$
By the choice of $u$ (not a cubic polynomial),
we have the following:
\begin{align}
&r_E(\hat u)=\xi \big{\vert}_E \ne 0 \quad 
\mbox{in} \quad H^0(\mathcal O_E(2)), \quad \mbox{and} \label{non-zero}\\
&\div (\xi \big{\vert}_E)=Z. \label{zero}
\end{align}
These respectively follow from 
the explicit description of $r_E$ 
in Claim \ref{claim:explicit description of r_E}
and the direct diagram chasing.

Before we start the computation, we observe one sheaf
inclusion $\mathcal O_C(2Z) \subset {\mathcal N_C}^{\vee}(3)(2Z)$.
We get the inclusion by taking the dual of the exact sequence 
of normal bundles
\begin{equation}\label{ses:normal bundles sequence}
0 \longrightarrow
\underbrace{\mathcal N_{C/S}}_{\cong \, \omega_C(1)} \longrightarrow
\mathcal N_C \longrightarrow
\underbrace{\mathcal N_S \otimes \mathcal O_C}_{\cong \, \mathcal O_C(3)} 
\longrightarrow 0
\end{equation}
and then tensoring with $\mathcal O_C(3)(2Z)$.
We see that the inclusion induces an injection between their $H^1$.
For the injectivity, it is enough to show that 
$\mathcal N_{C/S}^{\vee}(3)(2Z)
\cong \mathcal O_S(3\mathbf h+2E-C)\big{\vert}_C$
does not have global sections.
Indeed, we have
$$
(3\mathbf h+2E-C)\cdot C
=-(C-3\mathbf h-E)^2 -3\mathbf h \cdot (C-3\mathbf h-E) +2<0,
$$
since $C-3\mathbf h-E$ is nef (hence effective) and big.
Therefore we get the injection.

\begin{lem}\label{lem:obstruction}
Let $\varphi$, $u$, and $\xi$ be as above. 
Let $\mathbf t$ be the image of $\bar \delta(u)$ by the map
$H^1(\mathcal I_C(3)) \rightarrow H^1({\mathcal N_C}^{\vee}(3)(2Z))$.
Then we have the following:
\begin{enumerate}
 \item $\mathbf t$ is contained in 
       $H^1(\mathcal O_C(2Z)) \subset H^1({\mathcal N_C}^{\vee}(3)(2Z))$.
       Moreover, the cup product by $\cup_4$
       corresponding to $\varphi$ equals
       the cup product $\mathbf t \cup u$ by
       $$
       \cup_5: H^1(\mathcal O_C(2Z)) \times H^0(\mathcal O_C(3))
       \overset{\cup}{\longrightarrow} H^1(\mathcal O_C(3)(2Z)).
       $$
 \item Let $p:C \rightarrow \mathbb P^1$ be the projection from $Z$,
       and let $x',y'$ be two linearly independent global sections of
       $p^{\ast} \mathcal O_{\mathbb P^1}(1)=\mathcal O_C(\mathbf h -Z)$
       corresponding to $x,y$.
       Then $\mathbf t$ is represented by a 1-cocycle
       $$
       \dfrac{\xi}{xy} \in C^1(\mathfrak U_1,\mathcal O_C(2Z))
       =\Gamma(D(x')\cap D(y'),\mathcal O_C(2Z))
       $$
       with respect to the open affine covering
       $\mathfrak U_1=\left\{ D(x'),D(y') \right\}$ of $C$.
\end{enumerate}
\end{lem}
\noindent
{\sc Proof. \quad}
We compute the coboundary $\bar \delta(u)$ in $H^1(\mathcal I_C(3))$.
We recall the cubic equation $f=Ax+By$ defining $S$ (cf. \S 2.3).
By the smoothness of $S$, 
$$
 \mathfrak U_2: = \{ D(x),D(y),D(A),D(B) \}
$$
is an open affine covering of $\mathbb P^3$.
We compute $\bar \delta(u)$ 
by the {\v C}ech cohomology with respect to $\mathfrak U_2$.
By \eqref{presentation of $u$} $u$ is represented by
$\eta_1/x$ over $D(x)$ and $\eta_2/y$ over $D(y)$, where
$\eta_1,\eta_2$ are quartic polynomials such that 
$x\eta_2-y\eta_1=\xi f$.
Therefore, $\bar \delta(u)$ in $H^1(\mathcal I_C(3))$
is represented by 
$$
\bar \delta(u)=\dfrac{\eta_2}y - \dfrac{\eta_1}x 
=\dfrac{x \eta_2 -y \eta_1}{xy}
=\dfrac{\xi}{xy} f
$$
over $D(x) \cap D(y)$.
Thus $\bar \delta(u)$ is contained in the subsheaf
$\mathcal I_S(3) \subset \mathcal I_C(3)$
over $D(x) \cap D(y)$.
Restricting it to $C$, we see that $\mathbf t$ is contained in the
subsheaf $\mathcal O_C \subset {\mathcal N_C}^{\vee}(3)$
over $D(x) \cap D(y)$ and represented by $\xi/xy$ there.

On the other hand, the subcovering $\left\{ D(x),D(y) \right\}$ 
of $\mathfrak U_2$ covers whole $C$ except for $Z$.
Indeed, the two linear forms $\left\{x,y \right\}$ is a basis of
the pencil $P:=H^0(\mathcal O_C(1)(-Z))$ and
the fixed part of $P$ is exactly $Z$.
Therefore $D(x)=D(x')\setminus Z$, $D(y)=D(y')\setminus Z$, and
$\xi/xy$ gives a section of $\mathcal O_C(2Z)$ over $D(x')\cap D(y')$.
Now we make a change on the coverings of $C$.
We consider another open affine covering
$$
\mathfrak U_3:=\{ D(x'),D(y'),D(A), D(B) \}
$$
of $C$. Then both $\mathfrak U_1$ and $\mathfrak U_2$ are 
refinements of $\mathfrak U_3$.
There are isomorphisms between all {\v C}ech cohomology groups
$$
\check H^1(\mathfrak U_i,\mathcal O_C(2Z))
\quad (1\le i\le 3)
$$
induced by natural maps
$$
C^{\bullet}(\mathfrak U_1,\mathcal O_C(2Z))
\longleftarrow
C^{\bullet}(\mathfrak U_3,\mathcal O_C(2Z))
\longrightarrow
C^{\bullet}(\mathfrak U_2,\mathcal O_C(2Z))
$$
of {\v C}ech complexes with respect to $\mathfrak U_i$ ($1\le i\le 3$).
Moreover, by the above computation, 
we see that the 1-cocycle representing $\mathbf t$
can be taken from the one in  $C^1(\mathfrak U_3,\mathcal O_C(2Z))$,
and mapped to $\xi/xy$ in $C^1(\mathfrak U_1,\mathcal O_C(2Z))$.
Hence we have proved (2) and $\mathbf t \in H^1(\mathcal O_C(2Z))$.
Finally, we prove (1).
By the definition of $u$, the restriction of 
$\varphi$ to $\mathcal O_C \subset {\mathcal N_C}^{\vee}(3)$
as a homomorphism 
${\mathcal N_C}^{\vee}(3) \rightarrow \mathcal O_C(3)$ is 
the multiplication map by $u$.
The desired cup product is $\mathbf t \cup u$ by $\cup_5$.
\qed

\subsection{}
In this subsection, we show that the cup product $\mathbf t \cup u$ 
obtained in Lemma \ref{lem:obstruction} is non-zero.
For this purpose, first we show $\mathbf t\ne 0$.

\begin{claim}\label{claim:non-zero of 1st elt}
 $\mathbf t\ne 0$ in $H^1(\mathcal O_C(2Z))$.
\end{claim}
\noindent
{\sc Proof. \quad}
Since $\left\{ x,y \right\}$ is a basis of $H^0(\mathcal O_S(1)(-E))$,
by the base point free pencil trick, there exists an exact sequence 
$$
  \begin{CD}
   \mathbb K:
   0     @>>>
   \mathcal O_S(2E)  
   @>{\left[\begin{smallmatrix} y & -x \\
	    \end{smallmatrix}\right]}>>
   {\mathcal O}_S(1)(E)^{\oplus 2} 
   @>{\left[\begin{smallmatrix} x \\ y \\
	    \end{smallmatrix}\right]}>>
   \mathcal O_S(2)
   @>>> 0 \\
  \end{CD}
$$
of Koszul type. 
The restriction of $\mathbb K$ to $C$ is the exact sequence
$$
 \begin{CD}
  \mathbb K_C: 
  0     @>>>
  \mathcal O_C(2Z)  
  @>{\left[\begin{smallmatrix} y & -x \\
	   \end{smallmatrix}\right]}>>
  {\mathcal O}_C(1)(Z)^{\oplus 2} 
  @>{\left[\begin{smallmatrix} x \\ y \\
	   \end{smallmatrix}\right]}>>
  \mathcal O_C(2)
  @>>> 0.
 \end{CD}
$$
The restriction map $\mathbb K \rightarrow \mathbb K_C$ induces
$$
\begin{array}{ccccc}
 H^0({\mathcal O}_S(1)(E)^{\oplus 2}) &
  \mapright{\sigma} &
  H^0(\mathcal O_S(2)) &
  \mapright{\delta''} &
  H^1(\mathcal O_S(2E)) \\
 \mapdown{} && \mapdown{} && \mapdown{res} \\
 H^0({\mathcal O}_C(1)(Z)^{\oplus 2}) &
  \mapright{\gamma} &
  H^0(\mathcal O_C(2)) &
  \mapright{\delta'} &
  H^1(\mathcal O_C(2Z)). \\
\end{array}
$$
By the definition of the {\v C}ech coboundary map and
the description of $\mathbf t$ 
obtained in Lemma \ref{lem:obstruction} (2),
we have $\delta'(\xi\big{\vert}_C)=\mathbf t$.
Put $\hat {\mathbf t}:=\delta''(\xi \big{\vert}_S)$. 
Then $\hat {\mathbf t}$ is an element of $H^1(\mathcal O_S(2E))$
such that 
$\hat {\mathbf t}=\xi /xy$ over $D(x') \cap D(y')$ and 
$\hat {\mathbf t} \big{\vert}_C=\mathbf t$.

On the other hand, we obtain the exact sequence
$$
  \begin{CD}
   \mathbb K_E:
  0     @>>>
  \mathcal O_E(-2)  
  @>{\left[\begin{smallmatrix} A & B \\
	   \end{smallmatrix}\right]}>>
  {\mathcal O_E}^{\oplus 2} 
  @>{\left[\begin{smallmatrix} -B \\ A \\
	   \end{smallmatrix}\right]}>>
  \mathcal O_E(2)
  @>>> 0 \\
  \end{CD}
$$
as the restriction of $\mathbb K$ to $E$.
Here $A,B$ denote the quadratic polynomials
in the equation $f=Ax+By$ of $S$.
The restriction map $\mathbb K \rightarrow \mathbb K_E$ induces
$$
\begin{array}{ccccc}
   H^0({\mathcal O_S(1)(E)}^{\oplus 2}) &
   \mapright{\sigma} &
   H^0(\mathcal O_S(2)) &
   \mapright{\delta''} &
   H^1(\mathcal O_S(2E)) \\
   \mapdown{} && \mapdown{} && \mapdown{} \\
   H^0({\mathcal O_E}^{\oplus 2}) &
   \mapright{\varepsilon} &
   H^0(\mathcal O_E(2)) &
   \mapright{} &
   H^1(\mathcal O_E(-2)). \\
\end{array}
$$
It follows from $H^1(\mathcal O_S(2E-C))=0$ that 
the restriction map 
$H^1(\mathcal O_S(2E))\overset{res}{\rightarrow} H^1(\mathcal O_C(2Z))$
is injective. 
Hence it suffices to prove $\hat {\mathbf t}\ne 0$ for the claim.
Suppose that $\hat {\mathbf t}=0$ for contradiction.
Then $\xi\big{\vert}_S \in \im \sigma$ and hence
$\xi\big{\vert}_E \in \im \varepsilon$. This implies that
$\xi\big{\vert}_E$ is a linear combination of 
$A\big{\vert}_E$ and $B\big{\vert}_E$.
When we consider the divisors of zeros corresponding to 
$\xi\big{\vert}_E$ and $\langle A\big{\vert}_E,B\big{\vert}_E \rangle$,
this means $Z=\div(\xi\big{\vert}_E)$ (by \eqref{zero})
belongs to the restriction $\Lambda\big{\vert}_E$
of the conic pencil $\Lambda=|\mathbf h -E|$ to $E$
(cf. Remark \ref{rem:conic pencil}).
This contradicts Lemma \ref{lem:general $Z$}.
Thus $\hat {\mathbf t}\ne 0$.
\qed

\vskip 3mm

We next prepare an effective divisor $\Delta$ on $S$
which fills a gap between Mumford's case ($C \sim 4\mathbf h +2E$)
and our general case.
Consider the linear system 
$|C-4\mathbf h|$ ($\ne \emptyset$ by assumption) on $S$.
Since $(C-4\mathbf h-mE)\cdot E=-2+m<0$ if and only if $m < 2$,
$\Bs |C-4\mathbf h|$ contains $E$ with multiplicity two.
We take a member $\Delta$\footnote{
When $\Delta=0$, then $C \sim 4\mathbf h + 2E$ on $S$.
This is exactly the case of Mumford's example (\cite{Mumford}).
Taking $\Delta=0$ in our proof,
we have a proof for his case.
Thus Proposition \ref{prop:core proposition}
is a natural generalization of his example.}
of $|C-4\mathbf h-2E|$
which is disjoint from $E$ and fix it.
Then there exists a cup product map
$$
\cup_6: H^1(\mathcal O_C(2Z)) \times H^0(\mathcal O_C(3)(\Delta))
\longrightarrow H^1(\mathcal O_C(3)(2Z+\Delta)),
$$
which is compatible with $\cup_5$ via natural maps.
The last sheaf $\mathcal O_C(3)(2Z+\Delta)$ is isomorphic to 
the canonical line bundle $\mathcal O_C(K_C)$ by 
$$
\mathcal O_C(3)(2Z+\Delta)
\cong \mathcal O_S(3\mathbf h+2E+\Delta)\big{\vert}_C
\cong \mathcal O_S(C-\mathbf h)\big{\vert}_C 
\cong \mathcal O_C(K_C).
$$
Thus $\cup_6$ is the Serre duality cup pairing for $\mathcal O_C(2Z)$.
Now we consider an exact sequence
\begin{equation}\label{ses:canonical bundles}
 0 \longrightarrow 
  \underbrace{\mathcal O_S(-\mathbf h)}_{\cong \, \mathcal O_S(K_S)}
  \longrightarrow 
  \mathcal O_S(C-\mathbf h) 
  \longrightarrow 
  \underbrace{\mathcal O_S(C-\mathbf h)\big{\vert}_C}%
  _{\cong \, \mathcal O_C(K_C)}
  \longrightarrow 0. 
\end{equation}
Then the cup product map with its extension class $\mathbf e$,
which is the coboundary map of \eqref{ses:canonical bundles},
induces the next commutative diagram:
$$
\begin{CD}
 H^1(\mathcal O_C(2Z)) @. \times  @. H^0(\mathcal O_C(3)(\Delta)))
  @>{\cup_6}>> H^1(\mathcal O_C(3)(2Z+\Delta)) @. \cong H^1(K_C)\\
 @A{\res}AA @.  @VV{\cup \: \mathbf e}V @VV{\cup \: \mathbf e}V @. \\
 H^1(\mathcal O_S(2E)) @. \times  @. H^1(\mathcal O_S(3)(\Delta-C)))
  @>{\cup_7}>> H^2(\mathcal O_S(3)(2E+\Delta-C)) @. \cong H^2(K_S),\\
\end{CD}
$$
where $\res$ is the restriction map 
in the proof of Claim \ref{claim:non-zero of 1st elt}.
The last cup product $\cup_7$ is the Serre duality cup pairing 
for $\mathcal O_S(2E)$.

We have already got the non-zero element
$\hat {\mathbf t}$ of $H^1(\mathcal O_S(2E))$
such that $\hat {\mathbf t} \big{\vert}_C=\mathbf t$
in the proof of Claim \ref{claim:non-zero of 1st elt}.
By the commutativity of the diagram, we have
$(\mathbf t \cup u) \cup \mathbf e 
= \hat {\mathbf t}\cup (u \cup \mathbf e)$.
Since $H^1(\mathcal O_S(2E))$ is of dimension one,
by the Serre duality,
we have only to show that $u \cup \mathbf e \ne 0$ in 
$H^1(\mathcal O_S(3)(\Delta-C)))$ 
instead of $\mathbf t \cup u \ne 0$ in $H^1(K_C)$.

\begin{claim}
 $u \cup \mathbf e \ne 0$ in $H^1(\mathcal O_S(3)(\Delta-C))$.
\end{claim}
\noindent
{\sc Proof. \quad}
Suppose that $u \cup \mathbf e=0$ for contradiction.
Since the cup product map with $\mathbf e$ is
the coboundary map of the exact sequence
$$
0 \longrightarrow \mathcal O_S(3)(\Delta -C)
\longrightarrow \mathcal O_S(3)(\Delta) 
\longrightarrow \mathcal O_C(3)(\Delta)
\longrightarrow 0,
$$
there exists $\hat u'$ in $H^0(\mathcal O_S(3)(\Delta))$
such that $\hat u'\big{\vert}_C=u$.
Since $\Delta$ and $E$ are disjoint,
the image of $\hat u'$ by the restriction map
$$
H^0(\mathcal O_S(3)(\Delta+E))
\overset{r_E}{\longrightarrow} H^0(\mathcal O_E(2))
$$
is zero. 
Now we recall that $u$ has a lift $\hat u$ in $H^0(\mathcal O_S(3)(E))$
such that $r_E(\hat u)\ne 0$ by \eqref{non-zero}. 
Since $H^0(\mathcal O_S(3)(\Delta+E-C))\cong H^0(S,-\mathbf h -E)=0$,
we deduce $\hat u'=\hat u$ from 
$\hat u'\big{\vert}_C=\hat u\big{\vert}_C=u$.
This is a contradiction.
\qed

Therefore we complete 
the proof of Proposition \ref{prop:core proposition}.
\qed

\subsection{}
In this subsection, we give a technical remark to
Proposition \ref{prop:core proposition}.
In the proof of this proposition, 
the assumption that $C$ is a general member
of $|\mathbf D|$ was used only to prove Lemma \ref{lem:general $Z$}.
We characterize the members $C$ that do not satisfy
$Z=C\cap E \notin \Lambda\big{\vert}_E$,
where $\Lambda$ is the conic pencil $|\mathbf h -E|$ on $S$.
\begin{prop}
Let $C$, $E$, $Z$, and $\Lambda$ be as above.
Then the following two conditions are equivalent:
 $(1)$ $Z \not\in \Lambda\big{\vert}_E$; 
 $(2)$ $H^0(\mathcal O_C(1)(-2Z))=0$.
\end{prop}
\noindent
{\sc Proof. \quad}
Let us consider the commutative diagram of restriction maps:
 $$
 \begin{array}{ccc}
  H^0(\mathcal O_S(1)(-E)) & \mapright{r_1} & H^0(\mathcal O_E(2)) \\
  \mapdown{v_1} && \mapdown{v_2} \\
  H^0(\mathcal O_C(1)(-Z)) & \mapright{r_2} & H^0(\mathcal O_Z).
 \end{array}
 $$
The first condition is equivalent to the injectivity of the
the composite $v_2 \circ r_1$.
On the other hand, the second condition is equivalent to
the injectivity of $r_2$.
Therefore, it suffices to show that $v_1$ is an isomorphism.
In fact, we can easily check that $C -\mathbf h +E$ is nef and big.
This implies that $H^i(\mathcal O_S(1)(-E-C))=0$ for $i=0,1$.
Thus we have the equivalence.
\qed

\vskip 3mm

Suppose that $H^0(\mathcal O_C(1)(-2Z))\ne 0$.
Then there exists a plane $H$ 
which is tangential to $C$ at $Z$.
Let $Z=p+q$ where $p,q \in C$.
Then the tangents to $C$ at $p$ and $q$ 
are coplanar. (See Figure 1.)
\begin{figure}[h]\label{fig:tangents}
 \begin{center}
  \includegraphics[clip,scale=0.7]{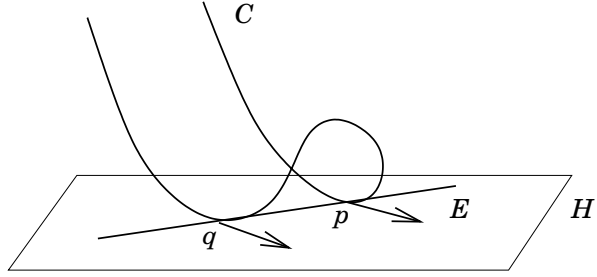}
  \caption{two tangents on a plane}  
 \end{center}
\end{figure}

When $Z \in \Lambda\big{\vert}_E$,
what can we say about the obstruction?
Let $C$ be such a special member of $|\mathbf D|$.
Then the reverse diagram chase in the proof of
Claim \ref{claim:non-zero of 1st elt}
shows $\mathbf t=0$.
Thus the cup product $\mathbf t \cup u$ by $\cup_5$ is zero.
Since $H^1(\mathcal N_C),H^1(\mathcal O_C(3))$ and 
$H^1(\mathcal O_C(3)(2Z))$ are all isomorphic via natural maps,
we deduce all the previous cup products by $\cup_i$ ($i \le 4$)
are zero.
Hence $\varphi \in H^0(\mathcal N_C)$ corresponding to 
$\mathbf t=0$ is not obstructed at the second order.
However, we will later see that $C$ corresponds to
a non-reduced point of the Hilbert scheme
(cf. Proposition \ref{prop:non-reduced components}).
This implies that  $\varphi$ is obstructed at the 
$n$-th order for some $n \ge 3$.

\section{An application to non-reduced components of the Hilbert scheme}

In this section, we apply Proposition \ref{prop:core proposition} 
to a problem on non-reduced components of the Hilbert scheme,
and prove the main theorem.
The theorem shows that a special case of Conjecture 
\ref{conj:Kleppe-Ellia} of Kleppe and Ellia is true.

Let $W$ be an irreducible closed subset of the Hilbert scheme
$H_{d,g}^S$ with $d\ge 3$.
Suppose that $W$ is maximal among all the irreducible closed subsets
of $H_{d,g}^S$ whose general member is contained in a 
smooth cubic surface. Let $C$ be a general member of $W$ 
and let $S$ be a general cubic surface containing $C$. 
Then we obtain a $7$-tuple $(a;b_1,\ldots,b_6)$ of integers satisfying
\begin{equation}\label{form:deg and gen}
 \left\{
  \begin{aligned}
   & a > b_1 \ge b_2 \ge \ldots \ge b_6 \ge 0, 
   \quad a\ge b_1+b_2+b_3, \\
   & d = 3a - \sum_{i=1}^6 b_i , \quad \mbox{and} \quad
   g = \binom{a-1}{2} - \sum_{i=1}^6
   \binom{b_i}{2}
  \end{aligned}
  \right.
\end{equation}
as the $\mathbb E$-multidegree of $C$. (See \S 2.2 for more detail.)

Conversely, suppose that a $7$-tuple $(a;b_1,\ldots,b_6)$ satisfying
\eqref{form:deg and gen} is given. 
If $\mathcal L$ is an invertible sheaf of this multidegree
on a smooth cubic surface $S$, then every general member of 
$|\mathcal L|$ is a smooth connected curve
by the conditions $a>b_1$ and $b_6 \ge 0$. 
Thus we have a non-empty irreducible closed subset $W$ of $H^S_{d,g}$ by
\begin{equation}\label{eqn:description of W}
 W := \left\{
     \mbox{$C \in H^S_{d,g} \ \big{\vert}$ $C \subset S$: a smooth cubic, } 
     \mathcal O_S(C) \cong \mathcal O_S(a;b_1,\ldots,b_6)
     \in \Pic S \right\}^{-},
\end{equation}
where ${}^{-}$ denotes the closure in $(H_{d,g}^S)_{\red}$.

\begin{defn}
 For a $7$-tuple $(a;b_1,\ldots,b_6)$ of integers satisfying
 \eqref{form:deg and gen}, we denote the above subset $W$ 
 of $H^S_{d,g}$ by $W_{(a;b_1,\ldots,b_6)}$.
\end{defn}

When $d > 9$, any general member $C$ of $W$ is contained in 
the unique cubic surface $S$, and furthermore, the above construction gives
one-to-one correspondence
$(a;b_1,\ldots,b_6) \leftrightarrow W_{(a;b_1,\ldots,b_6)}$ between 
the 7-tuples satisfying \eqref{form:deg and gen} and 
the maximal irreducible closed subsets $W$ of $H_{d,g}^S$
whose general member is contained in a smooth cubic surface
(cf. \cite[Remark 2]{Kleppe85}).
Thus to determine all irreducible components of $H_{d,g}^S$ 
whose general member is contained in a smooth cubic,
it suffices to solve the next problem:
\begin{prob}\label{prob:class'n problem}
Determine all $W_{(a;b_1,\ldots,b_6)}$ that are
irreducible components of $(H_{d,g}^S)_{\red}$.
\end{prob}
The above problem makes sense only when $g \ge 3d-18$.
This is because, as is found in \cite{Kleppe85},
$\dim W_{(a;b_1,\ldots,b_6)}=d+g+18$ when $d>9$, while
every irreducible component of $H_{d,g}^S$ is of dimension at least 
$4d$ ($=\chi(\mathcal N_C)$) from a general theory.
In what follows, we consider the above problem in the range
$$\Omega:=\{(d,g)\in \mathbb Z^2|d> 9, \ g\ge 3d-18\}.$$

Let $(d,g) \in \Omega$, let $W=W_{(a;b_1,\ldots,b_6)}$,
and let $C$ be a general member of $W$.
Then, we have natural inequalities
\begin{equation}\label{ineq:basic inequalities}
\dim W \le \dim_{[C]} H_{d,g}^S \le h^0(\mathcal N_C).
\end{equation}
If $\dim W = \dim_{[C]} H_{d,g}^S$, 
then $W$ is an irreducible component of $(H_{d,g}^S)_{\red}$.
$H_{d,g}^S$ is smooth at $[C]$ if and only if 
$\dim_{[C]} H_{d,g}^S = h^0(\mathcal N_C)$.
The exact sequence \eqref{ses:normal bundles sequence} induces 
$H^1(\mathcal N_C) \cong H^1(\mathcal O_C(3))$
because we have $H^1(\omega_C(1))=0$.
Therefore we get
\begin{equation}\label{eqn:imbed codim}
h^0(\mathcal N_C) - \dim W
= (4d + h^1(\mathcal O_C(3))) - (d+g+18)
=h^1(\mathcal I_C(3)). 
\end{equation}
Here the last equality follows from the exact sequence 
\eqref{les:cubic normality}.
By the same equality, in our case where $\dim W \ge 4d$, 
we always have 
\begin{equation}\label{ineq:h^1}
h^1(\mathcal O_C(3)) \ge h^1(\mathcal I_C(3)). 
\end{equation}

Let $S$ be the cubic surface containing $C$ and
let $\mathbf h$ be the class of hyperplane sections of $S$.
Then, as we saw in \S 2, the dimension $h^1(\mathcal I_C(3))$
can be computed from the fixed part $F$ of the linear system
$\Lambda_3:=|C-3\mathbf h|$ on $S$.
By the formula \eqref{form:fixed parts},
$F$ is empty (i.e. $\Lambda_3$ is free),
or a union of three kinds of (multiple) lines:
single, double, or triple.

\begin{lem}\label{lem:imbed codim}
 Let $(d,g)\in \Omega$, let $W=W_{(a;b_1,\ldots,b_6)} 
 \subset H_{d,g}^S$, and  let $C$ be as above.
 \begin{enumerate}
  \item If $d<12$, then $H^1(\mathcal I_C(3))=0$.
  \item If $d \ge 12$, then we have
	$$
	h^1(\mathcal I_C(3))=
	\sharp\left\{i|b_i=2\right\} 
	+3 (\sharp\left\{i|b_i=1\right\})
	+6 (\sharp\left\{i|b_i=0\right\}),
	$$ where $\sharp$ denotes the cardinality of a set.
	In particular, $H^1(\mathcal I_C(3))=0$ if and only if $b_6 \ge 3$.
 \end{enumerate}
\end{lem}
\noindent
{\sc Proof. \quad}
Let $S$, $\mathbf h$, and $\Lambda_3$ be as above. By the Serre duality,
we have
\begin{equation}\label{isoms:obst space}
H^1(\mathcal O_C(3))^{\vee} 
\cong H^2(\mathcal O_S(3\mathbf h -C ))^{\vee}
\cong H^0(\mathcal O_S(C-4\mathbf h)).
\end{equation}
Suppose $d<12$. Then the last cohomology group vanishes 
because $(C-4\mathbf h)\cdot \mathbf h=d-12$.
This implies $H^1(\mathcal I_C(3))=0$ by \eqref{ineq:h^1}.
Thus we proved (1).
Suppose $d\ge 12$. Then by the Riemann-Roch theorem on $S$,
we have $\chi(C-3\mathbf h)=g-2d+9\ge d-9>0$,
while $H^2(C-3\mathbf h)\cong H^0(2\mathbf h -C)^{\vee} =0$.
Therefore $C-3\mathbf h$ is effective.
Similarly, we have $(C-3\mathbf h)^2=2g-5d+25 \ge d-11 >0$.
By applying Corollary \ref{cor:vanishing thm2}
to $\Lambda_3$, we get the conclusion.
\qed

\vskip 3mm

When $b_6 \ge 3$ (i.e. $\Lambda_3$ is free), the lemma shows 
$H^1(\mathcal I_C(3))=0$.
This implies $h^0(\mathcal N_C)=\dim W$ by \eqref{eqn:imbed codim}.
Thus the following is obvious.
\begin{prop}[Kleppe \cite{Kleppe85}]\label{prop:reduced components}
 Let $(d,g)\in \Omega$ and let $W=W_{(a;b_1,\ldots,b_6)} \subset H^S_{d,g}$.
 If $b_6\ge 3$, then $H_{d,g}^S$ is generically non-singular along $W$.
 Moreover, $W$ is an irreducible component of $H_{d,g}^S$.
\end{prop}

When $d \ge 12$ and $b_6 \le 2$ (i.e. $\Lambda_3$ is non-free),
we have $h^1(\mathcal I_C(3))\ne 0$ by Lemma \ref{lem:imbed codim}.
So there may be some irreducible component $V$ which strictly
contains $W$. In this case, Problem 
\ref{prob:class'n problem} becomes non-trivial.
However, as long as we study the case where $h^1(\mathcal I_C(3))=1$,
the dichotomy between (A) and (B) described in the introduction (cf. \S 1)
makes the situation simple.
Now we give a proof of Theorem \ref{thm:main}.

\vskip 3mm

\noindent
\underline{\bf Proof of Main Theorem} \quad
Let $W$ be as in the statement.
Then $W=W_{(a;b_1,\ldots,b_6)}$ for some
$7$-tuple $(a;b_1,\ldots,b_6)$ satisfying \eqref{form:deg and gen}.
Lemma \ref{lem:imbed codim} shows that 
we have $h^1(\mathcal I_C(3))=1$
if and only if $d \ge 12$, $b_6=2$ and $b_5\ge 3$.
Thus the proof of the theorem reduces to the next proposition
which is an application of Proposition \ref{prop:core proposition}.

\begin{prop}\label{prop:non-reduced components}
 Let $d\ge 12$, let $g\ge 3d-18$, and
 let $W=W_{(a;b_1,\ldots,b_6)} \subset H_{d,g}^S$.
 If $b_6=2$ and $b_5\ge 3$, then
 $H_{d,g}^S$ is generically singular along $W$.
 Moreover, $W$ is an irreducible component of $(H_{d,g}^S)_{\red}$.
 Hence $H_{d,g}^S$ is non-reduced along $W$. 
\end{prop}
\noindent
{\sc Proof. \quad}
We check that any general member $C$ of $W$ satisfies the two conditions
(i) and (ii) of Proposition \ref{prop:core proposition}. 
The condition (i) is clearly satisfied with $E=E_6$ 
because of Lemma \ref{lem:zariski decomp} (iii).
Since $h^1(\mathcal I_C(3))=1$, 
we have $H^1(\mathcal O_C(3))\ne 0$ by \eqref{ineq:h^1}.
Therefore, the condition (ii) follows from \eqref{isoms:obst space}.

Since $C$ has an obstructed deformation
by Proposition \ref{prop:core proposition},
$H_{d,g}^S$ is singular at $[C]$ and we have
$\dim_{[C]} H_{d,g}^S < h^0(\mathcal N_C)$.
Consequently, we have $\dim W = \dim_{[C]} H_{d,g}^S$ in
\eqref{ineq:basic inequalities} from $h^1(\mathcal I_C(3))=1$.
Hence $W$ is an irreducible component of $(H_{d,g}^S)_{\red}$.
Moreover, since $H_{d,g}^S$ is singular at any general point of $W$,
$H_{d,g}^S$ is non-reduced along $W$.
\qed

Therefore the proof of Theorem \ref{thm:main} is completed.
\qed

\vskip 3mm

We give some example of non-reduced components of the Hilbert scheme.
\begin{exmp}\label{exa:non-reduced components}
Let $\lambda \ge 0$ be an integer. Then the subsets
\begin{align*}
&W_{(\lambda +12;\lambda+3,3,3,3,3,2)} \subset
H^S_{d,4d-37} \quad (d =2\lambda+19) \quad \mbox{and}\\
&W_{(\lambda+12;\lambda+4,3,3,3,3,2)} \subset
H^S_{d,\frac72 d -27} \quad (d =2\lambda +18)
\end{align*}
are irreducible components of $(H_{\mathbb P^3}^S)_{\red}$.
Moreover, $H_{\mathbb P^3}^S$ is non-reduced along each of them.
\end{exmp}

Theorem \ref{thm:main} shows that the next conjecture is true
whenever $h^1(\mathcal I_C(3))=1$ without the assumption that
$H^1(\mathcal I_C(1))=0$.
In fact, $H^1(\mathcal I_C(1))=0$ follows from $h^1(\mathcal I_C(3))=1$.

\begin{conj}[Kleppe \cite{Kleppe85}, Ellia \cite{Ellia}]
\label{conj:Kleppe-Ellia}
 Let $(d,g) \in \Omega$ and let $W$ be an irreducible closed subset
 of $H_{d,g}^S$ whose general member $C$ is contained in a smooth 
 cubic surface.
 Suppose that $W$ is maximal among all such subsets.
 If $H^1(\mathcal I_C(3))\ne 0$ and $H^1(\mathcal I_C(1))=0$, 
 then $W$ is an irreducible component of $(H_{d,g}^S)_{\red}$
 of dimension $d+g+18$. Moreover, $H_{d,g}^S$ is non-reduced along $W$.
\end{conj}

\begin{rem}\label{rem:Kleppe-Ellia}
This was originally conjectured by Kleppe in \cite{Kleppe85}
without the assumption of linearly normality ($H^1(\mathcal I_C(1))=0$).
He proved that the conjecture is true in the following two ranges:
$g> 7+(d-2)^2/8$ for $d \ge 18$, 
$g > -1+(d^2-4)/8$ for $14 \le d \le 17$.
When $d < 14$, we have $H^1(\mathcal I_C(3))=0$
by e.g. Lemma \ref{lem:imbed codim} (1) or \cite[Corollary 17]{Kleppe85}.
Hence he considered the conjecture in the range $d \ge 14$.
Later, Ellia \cite{Ellia} proved the conjecture 
for the wider range that $g>G(d,5)$ for $d\geq 21$.
Here $G(d,5)$ denotes the maximal genus of
curves of degree $d$, not contained in a quartic surface.
$G(d,5)$ nearly equals $d^2/10$ for $d \gg 0$.
Moreover, he gave a counterexample for linearly non-normal curves, and
suggested restricting the conjecture to linearly normal curves.

After the original version of this paper was submitted, 
the author learned that Kleppe \cite{Kleppe96} had made
further progress in proving the conjecture:
his result consists of a proof of the conjecture
for part of the case $h^1(\mathcal I_C(3))=1$
and that for part of the case $h^1(\mathcal I_C(3))=3$,
but does not cover our result
(cf. Example \ref{exa:non-reduced components}).
The method of his proofs is different from ours 
(cf. Remark \ref{rem:methods}).
\end{rem}

\begin{rem}\label{rem:methods}
To prove Conjecture \ref{conj:Kleppe-Ellia} for a given
$W=W_{(a;b_1,\ldots,b_6)} \subset H_{d,g}^S$,
it suffices to prove that
$W$ is a component of $(H_{d,g}^S)_{\red}$
because $H_{d,g}^S$ is automatically non-reduced along $W$
by the assumption $H^1(\mathcal I_C(3))\ne 0$.
In \cite{Kleppe85},\cite{Ellia} and \cite{Kleppe96},
the authors proved that $W$ is a component of $(H_{d,g}^S)_{\red}$
by contradiction.
First they assumed that a general member $C$ of $W$ is a specialization 
of curves contained not in a cubic but in a surface of
degree greater than three.
Then they got a contradiction 
by using a dimension count of a certain family of curves 
on a quartic (\cite{Kleppe85}, \cite{Ellia}),
or using the fact that the dimension of cohomology groups 
can only increase under specialization by semicontinuity
(\cite{Kleppe96}).
\end{rem}

Finally we remark that
$W$ in Conjecture \ref{conj:Kleppe-Ellia}
is not an irreducible component of $(H_{d,g}^S)_{\red}$ 
provided that $h^1(\mathcal I_C(1))\ne 0$.
This fact is obtained from the following,
whose proof is essentially given by
\cite[Remark VI.6]{Ellia} and \cite[Remark 2.10]{Dolcetti-Pareschi}.
\begin{prop}
[(Ellia \cite{Ellia}, Dolcetti-Pareschi \cite{Dolcetti-Pareschi})]
Let $(d,g)\in \Omega$ and let 
$W=W_{(a;b_1,\ldots,b_6)} \subset H_{d,g}^S$.
Suppose that $b_6=0$.
Then $W$ is not an irreducible component of $(H_{d,g}^S)_{\red}$.
\end{prop}

\vskip 3mm

\section*{Appendix (Irreducible components of $H_{d,g}^S$\\
whose general member is contained in a smooth quadric)}

We can naturally consider the same problem as Problem
\ref{prob:class'n problem} 
for curves contained in a smooth quadric surface 
$Q \cong \mathbb P^1 \times \mathbb P^1$ with bidegree 
$(a,b) \in \Pic Q \cong \mathbb Z^2$.
This problem is easier than Problem \ref{prob:class'n problem}.
One of the reason for this is that we have $H^1(Q,D)=0$
for any effective divisor $D$ on $Q$.

Let $d>4$ and $g\ge0$ be two integers.
For a pair $(a,b)$ of non-negative integers
satisfying $a+b=d$, $(a-1)(b-1)=g$ and $a \ge b >0$,
we define an irreducible closed subset $W_{(a,b)}$ of $H_{d,g}^S$
as follows:
$$
W_{(a,b)} := \left\{
     \mbox{$C \in H^S_{d,g} \ \big{\vert}$ $C \subset Q$: a smooth quadric, } 
     \mathcal O_Q(C) \cong \mathcal O_Q(a,b) 
     \in \Pic Q \right\}^{-}.
$$
Then $W_{(a,b)}$ is an irreducible closed 
subset of $H_{d,g}^S$ whose general member is contained 
in a smooth quadric surface and maximal among all such subsets.
We can easily see that $\dim W_{(a,b)}=2d+g+8$.
The next proposition shows that
$W_{(a,b)}$ is an irreducible component of $H_{d,g}^S$
if and only if $g \ge 2d-8$.
\begin{prop}
 Let $d>4$ and $g\ge 0$ be two integers, 
 and let $W_{(a,b)} \subset H_{d,g}^S$. 
 Then $H_{d,g}^S$ is generically non-singular along $W_{(a,b)}$. 
 Moreover, if $g\ge 2d-8$, then $W_{(a,b)}$ is an irreducible 
 component of $H_{d,g}^S$. 
 Otherwise, $W_{(a,b)}$ is a subvariety of $H_{d,g}^S$ of 
 codimension $2d-8-g$.
\end{prop}
\noindent
{\sc Proof. \quad}
Let $C$ be a general member of $W_{(a,b)}$
which is contained in a smooth quadric surface $Q$, 
and let $\mathbf h$ be the class of hyperplane sections of $Q$.
Then the exact sequence 
$0 \rightarrow \mathcal I_Q(2) 
\rightarrow \mathcal I_C(2) 
\rightarrow \mathcal I_{C/Q}(2)
\rightarrow 0$ induces 
$H^i(\mathcal I_C(2))\cong H^i(\mathcal I_{C/Q}(2))$ for $i=1,2$.
Therefore, we obtain 
\begin{equation}\label{isoms:appendix}
H^i(\mathcal I_C(2))
\cong H^{2-i}(\mathcal O_Q(a-4,b-4))^{\vee} \qquad (i=1,2)
\end{equation}
by $H^i(\mathcal I_{C/Q}(2)) 
\cong H^i(\mathcal O_Q(2\mathbf h -C))
\cong H^{2-i}(\mathcal O_Q(C-4\mathbf h))^{\vee}$.

First we assume that $g\ge 2d-8$. Since $g-2d+8=(a-3)(b-3)$,
we have $a\ge b>3$ when $g> 2d-8$, 
and we have $a=3$ or $b=3$ when $g=2d-8$.
Thus it follows from \eqref{isoms:appendix} that $H^1(\mathcal I_C(2))=0$.
By \cite[Theorem 1 (a)]{Kleppe85}, 
$W_{(a,b)}$ is a reduced component of $H_{d,g}^S$.
Next we assume that $g<2d-8$. This implies $b<3$ 
and hence we have $H^1(\mathcal O_C(2))
\cong H^2(\mathcal I_C(2))=0$ by \eqref{isoms:appendix}.
By \cite[Theorem 1 (b)]{Kleppe85}, 
$H_{d,g}^S$ is generically non-singular along $W_{(a,b)}$,
and the codimension of $W_{(a,b)}$ in $H_{d,g}^S$ is equal to 
$h^1(\mathcal I_C(2))=2d-8-g$.

Thus we conclude that $H_{d,g}^S$ is generically non-singular 
along $W_{(a,b)}$.
\qed

\end{document}